\documentclass{article}
\usepackage{amsthm}
\usepackage{color}
\usepackage{hyperref}
\hypersetup{colorlinks=true, 
 linkcolor=cyan,
 urlcolor=cyan}

\newcommand{\ds}{\displaystyle}
\newcommand{\bs}{\bigskip}
\newcommand{\ms}{\medskip}
\newcommand{\M}{{\cal M}}

\newcommand{\be}{\begin{enumerate}}
\newcommand{\ee}{\end{enumerate}}
\newcommand{\bi}{\begin{itemize}}
\newcommand{\ei}{\end{itemize}}

\newcommand{\abs}[1]{{\left|#1\right|}}
\newcommand{\pa}[1]{{\left(#1\right)}}
\newcommand{\cro}[1]{{\left[#1\right]}}

\catcode`\é=\active
\defé{\'e}

\catcode`\è=\active
\defè{\`e}

\catcode`\ê=\active
\defê{\^e}

\catcode`\û=\active
\defû{\^u}

\catcode`\â=\active
\defâ{\^a}

\catcode`\à=\active
\defà{\`a}

\catcode`\ù=\active
\defù{\`u}

\catcode`\ç=\active
\defç{\c c}

\catcode`\î=\active
\defî{\^{\i}}

\catcode`\ï=\active
\defï{\"{\i}}

\catcode`\ô=\active
\defô{\^o}

\catcode`\é=\active
\defé{\'e}

\font\ineg=msam8
\newcommand{\ie}{\mathrel{\hbox{\ineg 6}}}
\newcommand{\se}{\mathrel{\hbox{\ineg >}}}

\font\matcinq=msbm5
\font\matsept=msbm7
\font\matdix=msbm10
\newfam\matfam \scriptscriptfont\matfam=\matcinq \scriptfont\matfam=\matsept \textfont\matfam=\matdix
\newcommand{\mat}{\fam\matfam}
\newcommand{\N}{{\mat N}}
\newcommand{\Z}{{\mat Z}}
\newcommand{\R}{{\mat R}}
\newcommand{\C}{{\mat C}}
\newcommand{\SM}{{\mat S}}

\newcommand{\bib}[2]{\hbox{\hbox to 9mm{[#1]:\hfill} \hfill \hbox to 144mm{\vtop{\hsize=144mm#2\vfill\ms}\hfill}}}

\begin{document}


\begin{center}
\LARGE A unified approach for summation formulae
\end{center}
\ms

\begin{center}
\large Jean-Christophe Feauveau%
\footnote{Jean-Christophe Feauveau,\\
Professeur en classes préparatoires au lycée Bellevue,\\
135, route de Narbonne BP. 44370, 31031 Toulouse Cedex 4, France,\\
\href{mailto:Jean-Christophe.Feauveau@ac-toulouse.fr}{Jean-Christophe.Feauveau@ac-toulouse.fr}
}

\end{center}

\begin{center}
\large April 10, 2016
\end{center}

\bs

\textsc {Abstract}.
Summation formulae are classical tools in analysis: Taylor-MacLaurin,  Euler-MacLaurin, Poisson, Voronoï, Circle formulae\ldots
\ms

We will show how, from a single equation - referred to as the mother-equation - it is possible to unify these formulae 
and many others within a common formalism. Indeed, these formulae are paired up: every summation formula is associated with an asymptotic expansion. For example, the Euler-MacLaurin's formula turns out to be the asymptotic expansion associated with the Poisson's formula, the Taylor-MacLaurin's formula being of course the expansion of the function initially considered.
\ms

As is generally the case in sciences, this unifying concept is also generating new results. Asymptotic expansions for Voronoï and Circle formulae are presented first, then we show how to develop a M\"obius-Poisson's summation formula with its Euler-M\"obius-Poisson asymptotic expansion.
\bs
\bs

\bs

Classification A.M.S. 2010: 30B10, 30D10, 30E20.
\bs


{\textsc {I - Motivation}}
\ms

Various reasons, but first of all experimentation, lead theoretical physicists to search for a frame as synthetic as possible to give an account of physical phenomena. Maxwell's laws or the unification of interactions between elementary particles are examples of this trend.
\ms

Mathematics is not often offered the chance to be concrete. Nevertheless, the unification of various mathematical concepts is often fecund and generative of knowledge.
To convince oneself of this, it is enough to think about concepts such as equivalence class (congruence, rank of a linear application\ldots), isomorphism or functor.
\ms

More modestly, the purpose of this article is to show how some usual summation formulae can be unified, and how it is possible to generate new ones using the proposed formalism.
We do not look for `optimal' results in terms of hypotheses: we demonstrate the existence of a common source for them all, even if each of these expressions comes from `natural' developments and from various `historical' contexts.
\ms

The mother-equation giving access those summation formulae is
$$\forall y > 0, \ \ H[F](y) = \frac{1}{2i\pi}\int_{C_\infty} \Gamma(-s) K(-s)F^{(s)}(0)(y)^{s}ds. \leqno{(0)}$$

In this formula, $F$ is the function given as argument to the summation formula and $F^{(z)}$ its derivative of order $z\in \C$, $K$ is a kernel function and $C_\infty$ is a path  (in fact a limit of paths) specific to each summation formula. Our strategy consists in interpreting this integral in two ways (one inverse Mellin transform and a calculation of residues) to bring forward the summation formula. The same formalism allows, in addition, to derive an asymptotic expansion for any summation formula.
\ms

Equation $(0)$ thus appears as a generator of summation formulae.
\bs

\begin{center}
 \begin{tabular}{|p{2cm}|c|c|p{3cm}|p{3cm}|}
    \hline
    Summation \mbox{expression} & Path & Kernel K(s) & Summation \mbox{formula's name} & Asymptotic  \mbox{expansion's name}  \\
    \hline
    \hline
    $F(t)$ 
    
    & $\ds \Re(s) = -\frac{1}{2}$ & $1$ & Identity & Taylor-MacLaurin\\
    \hline
    $\ds \sum_{n=0}^{+\infty} F(nt)$ \vspace{1mm} & $\ds \Re(s) = -\frac{3}{2}$ & $\zeta(s)$ & Poisson & Euler-MacLaurin\\
    \hline
    $\ds \sum_{n=0}^{+\infty} d_n F(nt)$ \vspace{1mm} & $\ds \Re(s) = -\frac{1}{2}$ & $\zeta(s)^2$ & Voronoï & Euler-Voronoï\\
    \hline
    $\ds \sum_{n=0}^{+\infty} r_n F(nt)$ \vspace{1mm} & $\ds \Re(s) = -\frac{1}{2}$ & $\ds \zeta(s)\sum_{k=0}^{+\infty} \frac{(-1)^k}{(2k+1)^s}$ & Circle & Euler-Circle \\
    \hline
    $\ds \sum_{n=0}^{+\infty} \frac{\mu_n}{n} F(\frac{t}{n})$ & $\ds \abs{s} \to +\infty$ & $\ds \frac{(2\pi)^{-s}}{\zeta(-s+1)}$ & M\"obius-Poisson & Euler-M\"obius-Poisson\\
    \hline
  \end{tabular}
\end{center}
\bs
Some names in the previous table are not listed; the proposed names are then those that seemed the most natural according to their meaning.
\ms

The Taylor-MacLaurin case is separate: it corresponds to the sequence of weights $ (1,0,0, \ldots)$ in the summation. This expression is not really a summation - it is simply the function $F$ itself - but will give rise to an asymptotique expansion: the Taylor-MacLaurin's formula.
\ms

For all other expressions in the table, we will establish in order
\bi
\item
an asymptotic expansion with a $N^{{th}}$-order remainder,
\item
a proper summation formula that links the formula under study to the Fourier transform of $F$.
\ei

For reasons linked to the definition of a function's fractional derivative, we will only consider, in our asymptotic expansions, some functions $F: \R_-\to \C$. The variable, always noted $t$, will therefore be non-positive.
\ms

The   summation formulae being discussed naturally act on the even part of the real functions. That is why, in this framework, all functions are assumed even, with a non-negative variable $y$. We switch from a variable to the other one by $t =-y$.
\ms

Finally, in order to understand the common mechanism underlying these equations, we assume a fairly broad hypothesis on the functions $F$, namely that they be $C ^\infty$ functions with rapidly decreasing  derivatives.
\ms
\bs

  
{\textsc{II - Fractional derivation}}
\ms

For $a\in \R$, we call $\SM(]-\infty,a],\C)$ the space of functions $F \in C^\infty(]-\infty,a],\C)$ satisfying $\ds \lim_{t\to -\infty} \abs{t}^m F^{(n)}(t) = 0$ for all integers $n$ et $m$.
\ms

For any fixed $F \in \SM(]-\infty,a],\C)$ and $t \in ]-\infty,a]$, we define:
$$\forall s\in\C, \ \Re(s) < 0, \ \ F^{(s)}(t) = \frac{1}{\Gamma(-s)}\int_{-\infty}^t (t-u)^{-s-1} F(u)du = \frac{1}{\Gamma(-s)}\int_0^{+\infty} u^{-s-1} F(t-u)du.\leqno{(1)}$$
Through a simple calculation, for $s\in-\N ^\ast$, we obtain the iterated primitives of $F$ (with no constants of integration in $-\infty$). We have a classical extension of the notion of integration when $\Re(s) < 0$.
\ms

A translation by $t$ allows us to bring back the study of the punctual properties of fractional derivatives around $t = 0$. Without loss of generality, we consider this case in the following.
\bs

\textbf{Theorem:} \textit{For any $F \in \SM(]-\infty,0],\C)$, the function $s\mapsto F^{(s)}(0)$ defined on $\Re(s) < 0$
can uniquely be extended to a holomorphic function on $\C$. The extended function coincides with the usual $s^{th}$ derivative of $F$ at $0$ when $s\in \N$.}
\bs

\begin{proof}[Proof]
To establish this result, let us define
$\ds D(F)(s) := \frac{1}{\Gamma(-s)}\int_0^{+\infty} u^{-s-1} F(-u)du$.\\

For $\Re(s) < 0$, we get for any $n\in\N$
$$D(F)(s) = \frac{1}{\Gamma(-s)}\int_0^{+\infty} u^{-s-1} F(-u)du = \frac{1}{\Gamma(n-s)}\int_{-\infty}^0 F^{(n)}(u) (-u)^{n-s-1}du = D(F^{(n)})(s-n),$$
which proves the holomorphic extension for $\Re(s) < n$ then on $\C$ by uniqueness. 

Finally, for $n\in\N$, \ $\ds D(F)(n) = D(F^{(n+1)})(-1) = \int_0^{+\infty} F^{(n+1)}(-u)du = F^{(n)}(0)$.\\

This justifies writing $D(F)(s) = F^{(s)}(0)$ for $s\in\C$, and $(F^{(n)})^{(s)}(0) = F^{(n+s)}(0)$ for $n\in\N$.
\end{proof}
\ms

\textit {Remark:}
Working with $t \ie 0$ may seem unnatural but is linked to the definition used for the fractional derivative: that of Riemann-Liouville. Weyl's definition would allow us to work on $\R_+$, but the associated operator would no longer be the derivative $D$ but $-D$ and we would lose the intuitive character in the notation $F^{(s)}(0)$. And this one turns out to be very useful (cf. the Voronoï and Circle formulae).
\ms

A detailed study on fractional derivation and its various definitions can be found in [Ross].
\bs


{\textsc{III - Mellin Transform}}
\ms

If $F$ is defined on $\R_+$, its Mellin transform is [Tit1]
$$\M[F(u)](s) = \ds \int_0^{+\infty} v^{s-1} F(v)dv.\leqno{(2)}$$
In general, this transformation is only defined with certainty for $s$ in a strip such as $\sigma_0 <\Re(s) < \sigma_1$. 
It is typically the case when assuming  $F(u) = o(u^{-\sigma_1})$ in $+\infty$ 
and $F(u) = o(u^{-\sigma_0})$ in $0^+$.
\ms

Under this hypothesis, if for $c \in ]\sigma_0,\sigma_1[$ \, the function $t\mapsto \M [Fr (u)] (c+it) $ 
is integrable over $\R$ and $F$ is continuous, then we have the inversion formula [Tit1]
$$\forall v\in \R_+^\star, \ \ F(v) = \frac{1}{2i\pi}\int_{\Re(s)=c}  \M[F(u)](s) v^{-s}ds.\leqno{(3)}$$
\textit {Notation:} integration over \ $\Re (s) = c$ \ is to be taken along the orientated vertical line$] c-i\infty, c+i\infty [$.
\ms

Thus, for $F\in \SM(]-\infty,0],\C)$, we have
$$\forall s\in\C, \ \Re(s) < 0, \ \ F^{(s)}(0) = \frac{1}{\Gamma(-s)}\M[F(-u)](-s)$$
and
$$\forall c > 0, \ \forall v \ie 0, \ \ F(v) = \frac{1}{2i\pi}\int_{\Re(s) = c}  \Gamma(s)F^{(-s)}(0) (-v)^{-s}ds.$$
\bs\ms


{\textsc{IV - A Ramanujan's formula}}
\ms

Before studying usual summation formulae, let us see how fractional derivation and Mellin transform formally lead to a known formula.
According to Hardy, Ramanujan was very fond of this formula and used it diligently [Edwa].
\ms

For $F\in \SM(]-\infty,0],\C)$ assumed to be equal to its Taylor expansion over $]-\infty,0]$, we set $\Phi(s) = \ds \frac{F^{(s)}(0)}{\Gamma(s+1)}$. A Taylor-MacLaurin expansion for $F(-x)$ in the definition of $F^{(-s)}(0)$, and Euler's reflection formula $\ds \Gamma(s)\Gamma(1-s) = \pi (\sin(\pi s))^{-1}$, together lead to
$$\frac{\pi}{\sin(\pi s)} \Phi(-s) = \int_0^{+\infty} x^{s-1}\sum_{n=0}^{+\infty} (-1)^n x^n \Phi(n)dx,$$
which is actually Ramanujan's formula.
\bs\ms


{\textsc{V - Taylor-MacLaurin's formula}}
\ms

Through this case - the simplest - of asymptotic expansion, we introduce elements coming into play later. The result is the following.
\ms

For $F\in \SM(]-\infty,0],\C)$, we have
$$\forall N\in\N, \ \ \forall t\in \R_-, \ \ F(t) = \sum_{n=0}^{N} \frac{F^{(n)}(0)}{n!}t^n + \frac{1}{2i\pi} 
\int_{\Re(s)=N+\frac{1}{2}} \Gamma(-s)F^{(s)}(0) (-t)^{s}ds. \leqno{(4)}$$
\ms

Indeed, for $(a,b)\in\R^2$, $a < b$ and $T\in\R_+$, we denote by $R_{a,b,T}$ the rectangular path connecting the points $[a+iT,a-iT,b-iT,b+iT]$ 
in this order.
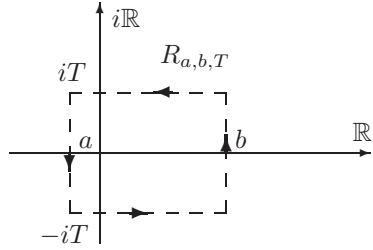
\ms

$$
\setlength{\unitlength}{0.8cm}
\begin{picture}(6,4)(-3,-2)

\put(3.2,0.2){$\R$} 
\put(-0.8,2.1){$i\R$}

\put(-1.35,0.09){$a$}
\put(1.25,0.09){$b$}

\put(-1.7,1.2){$iT$}
\put(-2,-1.5){$-iT$}

\put(0,1.5){$R_{a,b,T}$}

\put(-2.5,0){\vector(1,0){6}} 
\put(-1,-1.5){\vector(0,1){4}}

\multiput(-1.5,1)(0.4,0){7}{\line(1,0){0.2}} 
\multiput(-1.5,-1)(0.4,0){7}{\line(1,0){0.2}}

\multiput(-1.5,-1)(0,0.45){5}{\line(0,1){0.2}}
\multiput(1.1,-1)(0,0.45){5}{\line(0,1){0.2}}

\thicklines
\put(-1.5,0){\vector(0,-1){0.3}} 
\put(1.1,0){\vector(0,1){0.3}} 

\put(-0.5,-1){\vector(1,0){0.3}} 
\put(0.2,1.01){\vector(-1,0){0.3}} 

\end{picture}
$$

Cauchy's residue theorem gives, for any $T > 0$,
$$\frac{1}{2i\pi} \int_{R_{-\frac{1}{2},N+\frac{1}{2},T}} \Gamma(-s)F^{(s)}(0) (-t)^{s}ds = -\sum_{n=0}^{N} \frac{F^{(n)}(0)}{n!}t^n$$
since the only poles are due to $\Gamma(s)$, that they are located at $ \{1, \ldots, N \} $, and that their residues are
$\ds Res(\Gamma,-k) = (-1)^k / k!$     
as shown by the reflection formula.
\ms

It remains to control the integrals on horizontal segments of $R_{-\frac{1}{2}, N +\frac{1}{2}, T} $ and to interpret the integrals on vertical segments.
\ms

For $n\in\N$,
 $$\Gamma(-s) F^{(s)}(0) {\prod_{k = 0}^{n}(-s+k)} = \ds
 \int_{-\infty}^0 (-u)^{-s+n}F^{(n+1)}(u)du,$$
and since $F^{(n+1)}$ is  rapidly decreasing, this implies the existence of $M_n > 0$ such that $\abs{\Gamma(-s) F^{(s)}(0)} \ie \ds \frac{M_n}{\abs{\Im(s)}^n}$ for $\Re(s) \ie n+\frac{1}{2}$ and $\abs{\Im(s)} \se 1$.

The integer $N\in\N^*$ being fixed, we get $\ds \lim_{T\to +\infty} 
\int_{-1/2 + iT}^{N+1/2 + iT}
\Gamma(-s)F^{(s)}(0) (-t)^{s}ds = 0$ for $n=N$, and the same is true on the segment $[-\frac{1}{2} - iT,N+\frac{1}{2} - iT]$. 
This argument will be used again later.
\ms

The Mellin inversion formula then brings
$$\matrix{
\forall t\in \R_-, \ & F(t) & = & \ds \frac{1}{2i\pi} \int_{\Re(s)=-\frac{1}{2}} \Gamma(-s)F^{(s)}(0) (-t)^{s}ds\hfill\cr
& & = & \ds  \sum_{n=0}^{N} \frac{F^{(n)}(0)}{n!}t^n + \frac{1}{2i\pi} \int_{\Re(s)=N+\frac{1}{2}} \Gamma(-s)F^{(s)}(0) (-t)^{s}ds.\hfill\cr}$$
\bs

We obtain at the same time an expression for the $N$th-order remainder of Taylor-MacLaurin's formula:
$$\forall t \ie 0, \ \ Rt_N(F)(t) = \frac{1}{\Gamma(N+1)}\int_0^t F^{(N+1)}(u)(t-u)^n dt = \frac{1}{2i\pi} \int_{\Re(s) = c} \Gamma(-s)F^{(s)}(0) (-t)^{s}ds$$
for all $c\in ]N,N+1[$, for lack of pole on $N < \Re(s) < N+1$.
\ms

Thus, we derived the Taylor-MacLaurin formula from a calculation of $(0)$ in two different ways: a calculation of residues and an inversion of Mellin transform.\bs
\ms


{\textsc{VI - Euler-MacLaurin's formula}}
\ms

For $F\in \SM(]-\infty,0],\C)$ and  $t < 0$, the formula reads
$$\sum_{n=0}^{+\infty} F(nt) = -F^{(-1)}(0)t^{-1} + \frac{F(0)}{2} + \sum_{k=1}^{N} (-1)^k\frac{2\zeta(2k)}{(2\pi)^{2k}} F^{(2k-1)}(0) t^{2k-1} + Re_N(F).$$

The $\zeta$ functional equation is
$$\forall s\in\C-\{0,1\}, \ \ \ \ \xi(1-s) = \xi(s), \ \ {\rm with} \ \ \xi(s) = \pi^{-s/2}\Gamma(s/2)\zeta(s),\leqno{(5)}$$
and Euler-MacLaurin formula can be rewritten as
$$\matrix{
\ds \sum_{n=1}^{+\infty} F(nt) & = & \ds -F^{(-1)}(0)t^{-1} - \frac{F(0)}{2}  + \sum_{k=1}^{N} \frac{\zeta(1-2k)}{(2k-1)!} F^{(2k-1)}(0) t^{2k-1} + Re_N(F)\hfill\cr
& = &\ds -F^{(-1)}(0)t^{-1} + \sum_{k=0}^{2N-1} \frac{\zeta(-k)}{k!} F^{(k)}(0) t^{k} + Re_N(F).\hfill\cr}$$
\ms

We are going to obtain this classical result from a calculation of residue. Poles and zeros of $\Gamma$ lead us to study

$$H(t) = \frac{1}{2i\pi}\int_{R_{-\frac{3}{2},2N}} \Gamma(-s) \zeta(-s)F^{(s)}(0)(-t)^{s}ds$$

where $R_{-\frac{3}{2},2N}$ is the limit when $T$ tends to $+\infty$ of the rectangular paths $R_{-\frac{3}{2},2N,T}$ with vertices $[-\frac{3}{2}+iT,-\frac{3}{2}-iT,2N-iT,2N+iT]$.
The poles of the function being integrated are those of $\Gamma(-s)\zeta(-s)$. They are simple and their residues are:
\be
\item[{$\bullet$}]
$\ds Res(\Gamma(-s)\zeta(-s),2n-1) = \frac{\zeta(-2n+1)}{(2n-1)!}$ 
for $n\in\{1,\ldots,N\}$ ;
\item[{$\bullet$}]
 $\ds Res(\Gamma(-s)\zeta(-s),0) = \frac{1}{2}$ \, at $0$ ;
\item[{$\bullet$}]
$\ds Res(\Gamma(-s)\zeta(-s),1) = -1$ \, at $1$.
\ee

Hence, the calculation of residues gives us
$$H(t) = F^{(-1)}(0)t^{-1} - \sum_{k=0}^{2N-1} \frac{\zeta(-k)}{k!} F^{(k)}(0) t^{k}.$$

From $\ds \int_0^{+\infty} F(-nu)u^{s-1}du = \frac{1}{n^{s}}\int_0^{+\infty} F(-u)u^{s-1}du$, we deduce
$$\forall \Re(s) > 1, \ \ \zeta(s)\Gamma(s)F^{(-s)}(0) = \int_0^{+\infty} \sum_{n=1}^{+\infty} F(-nu)u^{s-1}du.\leqno{(6)}$$

At last, by inversion of Mellin transform, we obtain
$$\matrix{
 &\ds \sum_{n=1}^{+\infty} F(nt) & = & \ds \frac{1}{2i\pi} \int_{\Re(s)=\frac{3}{2}} \Gamma(s)\zeta(s)F^{(-s)}(0)(-t)^{-s}ds \hfill\cr
\ds (7) & & = & \ds \frac{1}{2i\pi} \int_{\Re(s)=-\frac{3}{2}} \Gamma(-s)\zeta(-s)F^{(s)}(0)(-t)^{s}ds \hfill\cr
(8) & & = & \ds -F^{(-1)}(0)t^{-1} + \sum_{k=0}^{2N-1} \frac{\zeta(-k)}{k!} F^{(k)}(0) t^{k} + \frac{1}{2i\pi} \int_{\Re(s)=2N} \Gamma(-s)\zeta(-s)F^{(s)}(0)(-t)^{s}ds. \hfill\cr}$$

The equation $(8)$ is justified by the integrals on the horizontal segments $[-\frac{3}{2} \pm iT,2N \pm iT]$ converging to 0 when $T$ tends towards $+\infty$.
\bs

Indeed, according to [Tit2], for any $\delta  < 0$, there exists $M_\delta > 0$ such as $\abs{\zeta(\sigma+ i t)} \ie M_\delta \abs{t}^{\frac{1}{2}-\delta}$ for $\sigma \se \delta$ and $\abs{t} \se 1$. The integer $N\in\N ^*$ being fixed, we choose here $\delta =-2N$ and in the argument of paragraph V, the choice $n=2N+1$ allows us to conclude:
there exists $M_n > 0$ such as $\abs{\Gamma(-s) F^{(s)}(0)} \ie \ds \frac{M_n}{\abs{\Im(s)}^n}$ for $\Re(s) \ie n+\frac{1}{2}$ and $\abs{\Im(s)} \se 1$.
\bs

Once more, an expression for the Euler-MacLaurin's order $N^{th}$-order remainder is obtained:
$$\forall t \ie 0, \ \ Re_N(F)(t) = \frac{1}{2i\pi} \int_{\Re(s)=c} \Gamma(-s)\zeta(-s)F^{(s)}(0) (-t)^{s}ds$$
for all $c\in ]2N-1,2N+1[$, assuming an absence of pole on $2N-1 < \Re(s) < 2N+1$.
Evaluating $H$ by a calculation of residues and its interpretation as Mellin transform therefore allows us to retrieve the Euler-MacLaurin formula.
\bs
\ms

Let us note the possibility to access a finite summation through the subtraction of two functions of $(8)$.

For $p< q$ two integers and $F\in \SM(]-\infty,0],\C)$, we can apply the previous expression to $F_{p} - F_{q}$, with $F_k(u) = F(u+kt)$, $t$ being fixed in $\R_-$. We find
$$\matrix{
\ds \sum_{n=p}^{q-1} F(nt) & = & \ds \frac{1}{2i\pi} \int_{\Re(s)=-\frac{3}{2}} \Gamma(-s)\zeta(-s)[F^{(s)}(pt) - F^{(s)}(qt)](-t)^{s}ds \hfill\cr
& = & \ds -(F^{(-1)}(pt)-F^{(-1)}(qt)) t^{-1} 
 + \sum_{k=1}^{N} (-1)^k\frac{2\zeta(2k)}{(2\pi)^{2k}} [F^{(2k-1)}(pt)-F^{(2k-1)}(qt)] t^{2k-1}\hfill\cr 
& & \ds {\hskip 1cm}
+ \frac{F(pt)- F(qt)}{2} 
+ \frac{1}{2i\pi} \int_{\Re(s)=2N} \Gamma(-s)\zeta(-s)[F^{(s)}(pt) - F^{(s)}(qt)](-t)^{s}ds. \hfill\cr\\}$$

Finally, for $t=-1$ et $G(u) = F(-u)$,
$$\frac{G(p)+G(q)}{2}+\sum_{n=p+1}^{q-1} G(n) = \int_p^q G(u)du + \sum_{k=1}^{N} (-1)^{k+1}\frac{2\zeta(2k)}{(2\pi)^{2k}} [G^{(2k-1)}(p)-G^{(2k-1)}(q)] + Re_{p,q}(G)$$
with $\ds Re_{p,q}(G) := \frac{1}{2i\pi} \int_{\Re(s)=c} \Gamma(-s)\zeta(-s)[F^{(s)}(-p)-F^{(s)}(-q)]ds$
for any $c\in ]2N-1,2N+1[$.
\ms

This is one of the usual forms of the Euler-MacLaurin's formula, up to the unusual form of the remainder.
\ms
\bs


{\textsc{VII - Poisson's formula}}
\ms

Under assumptions on $F$ to be specified later, the general form of Poisson's formula reads
$$\forall x\in\R, \ \ \forall y\in \R^\ast_+, \ \ \sum_{n\in\Z} F(ny+x) = \frac{1}{y} \sum_{n\in\Z} \widehat{F}\pa{\frac{2\pi n}{y}} e^{2i\pi nx/y},$$
which leads to
$$\forall y\in \R^\ast_+, \ \ \sum_{n\in\Z} F(ny) = \frac{1}{y} \sum_{n\in\Z} \widehat{F}\pa{\frac{2\pi n}{y}}. \leqno{(9)}$$

Those two equations are equivalent. The parameter $y$ is interpreted as the period of the function $F$.
\ms

The summation formula $(9)$ acts in fact only on the even part of $F$. By linearity, it is enough to prove it for $F$ even and real.
\ms

Let $F\in \SM (]-\infty, 0], \C) $ be extended into an even function on $\R$. 
It is easy to check using usual proofs that the  Poisson formula $(9)$ is proved in that case. 
It is now a matter of proving $(9)$ using the mother-equation $(0)$.
\ms

Following an approach found in [Mill], we are going to use the functional equation $(5)$ of $\zeta$.
Given that $F$ is even, the equality to be derived is
$$\forall y \in \R^\ast_+, \ \ \frac{F(0)}{2} + \sum_{n=1}^{+\infty} F(ny) = \frac{1}{y}\pa{\frac{\widehat{F}(0)}{2} + \sum_{n=1}^{+\infty} \widehat{F}\pa{\frac{2\pi n}{y}}}.$$

The starting point is expression $(7)$ used to obtain Euler-MacLaurin but with the $y > 0$  variable\\

$\matrix{
(10)   & \ds \sum_{n=1}^{+\infty} F(ny) & = &\ds \frac{1}{2i\pi}\int_{\Re(s) = -\frac{3}{2}} \Gamma(-s) \zeta(-s)F^{(s)}(0)y^{s}ds\hfill\cr
 & & = & \ds -\frac{1}{2}F(0) + \frac{1}{y}F^{(-1)}(0) + \frac{1}{2i\pi}\int_{\Re(s) = \frac{1}{2}} \Gamma(-s) \zeta(-s)F^{(s)}(0)y^{s}ds.\hfill\cr }  $\\

By a change of variable, it follows
$$\sum_{n=1}^{+\infty} F(ny) + \frac{1}{2}F(0) - \frac{1}{y}F^{(-1)}(0) = \frac{1}{2i\pi}\int_{\Re(s) = \frac{3}{2}} \Gamma(1-s) \zeta(1-s)F^{(s-1)}(0)y^{s-1}ds.$$

From \, $\ds \widehat{F}(x) = 2\int_0^{+\infty} F(u)\cos(xu)du$ \, 
we deduce \, 
$\ds \widehat{F}^{(s)}(0) = 2\cos(\pi s/2)\int_{0}^{+\infty} u^s F(u)du$,\, 
which can be rewritten as
$$\widehat{F}^{(-s)}(0) = 2\cos(\pi s/2)\Gamma(1-s)F^{(s-1)}(0) \leqno{(11)}$$
when $\Re(s) < 1$, then by analytical extension to $\C$ (for $s\in \N ^*$, $F ^ {(s-1)} (0) $ and $\cos (\pi s/2) $ alternately compensate the poles of $\Gamma$).
\ms

Moreover, Legendre duplication's formula $\ds \Gamma(s) = (2\pi)^{-1/2} 2^{s-1/2}\Gamma(s/2)\Gamma((s+1)/2)$ associated with the $\zeta$ functional equation leads to

$$\zeta(1-s) = 2(2\pi)^{-s}\Gamma(s)\cos(\pi s/2)\zeta(s).\leqno{(12)}$$

Given $(10)$, we find
$$\matrix{
\ds \sum_{n=1}^{+\infty} F(ny) + \frac{1}{2}F(0) - \frac{1}{y}F^{(-1)}(0) & = & \ds \frac{1}{2i\pi}\int_{\Re(s) = \frac{3}{2}} \Gamma(s)\zeta(s)\widehat{F}^{(-s)}(0)(2\pi)^{-s}y^{s-1}ds.\hfill\cr
& = &\ds \frac{1}{y}\sum_{n=1}^{+\infty} \widehat{F}(\frac{2\pi n}{y}).\hfill\cr}
$$
By remarking that $\frac{1}{t}F^{(-1)}(0) = \frac{1}{2t}\widehat{F}(0)$, formula $(7)$ is found.
Poisson's formula thus arise from a calculation of our mother-equation $(0)$ in two different ways: the Mellin's inversion formula and the Cauchy's residue theorem.
\bs\ms


{\textsc{VIII - Euler-Voronoï and Voronoï formulae}}
\ms

We note $d_n$ the number of divisors of $n$. Thus, for $\Re(s) > 1$, \,
 $\ds \sum_{n=1}^{+\infty} \frac{d_n}{n^s} = \zeta(s)^2$.

Under various hypothesis on function $F$ [Heja], [Endr], [Mill], we have the Voronoï formula in which $K_0$ and $Y_0$ denote the ordinary Bessel functions:

$$\sum_{n=1}^{+\infty} d_n F(n) = \int_0^{+\infty} (\ln(x)+2\gamma)F(x)dx + \frac{F(0)}{4} + 2\pi \sum_{n=1}^{+\infty} d_n\int_0^{+\infty}\cro{\frac{2}{\pi}K_0(4\pi\sqrt{nx}) - Y_0(4\pi\sqrt{nx})}F(x)dx. \leqno{(13)}$$
\ms

To derive this expression, following the example of the Euler-MacLaurin's formula, and for
$F\in \SM(]-\infty,0],\C)$ we introduce the function
$$\forall t < 0, \ \ H(t) = \frac{1}{2i\pi}\int_{R_{-\frac{3}{2},2N}} \Gamma(-s) \zeta(-s)^2F^{(s)}(0)(-t)^{s}ds \leqno{(14)}$$
where $R_{-\frac{3}{2},2N}$ is the limit of rectangular paths with vertices $[-\frac{3}{2}+iT,-\frac{3}{2}-iT,2N-iT,2N+iT]$.
\bs

\be
\item[{$\bullet$}]
At $2n-1$, $n\in-\N$, the pole is simple and $Res(\Gamma(-s)\zeta(-s)^2,2n-1) = \ds \frac{\zeta(-2n+1)}{(2n-1)!}$. 
\item[{$\bullet$}]
At $0$, the pole is simple and $Res(\Gamma(-s)\zeta(-s)^2,1) = \ds -\frac{1}{4}$. \\
\item[{$\bullet$}]
At $-1$, finally, there is a pole of order $2$:

$$\zeta^2(-s) = \frac{1}{(s+1)^2} - \frac{2\gamma}{s+1} + O(1), $$
$$ (-t)^{s} = -\frac{1}{t}(1+\ln(-t)(s+1) + O((s+1)^2)$$
and
$$F^{(s)}(0)\Gamma(-s) = F^{-1}(0) - (s+1)\int_0^{+\infty}F(-u)\ln(u)du + O((s+1)^2).$$
Hence, \, $\ds Res(\Gamma(-s)\zeta(-s)^2F^{(s)}(0)(-t)^s,1) = -\int_0^{+\infty}F(-u)\ln(\frac{u}{-t})du - 2\gamma\int_0^{+\infty}F(-u)du$.
\ee

The argument leading to $(6)$ and $(8)$ reads here
$$\matrix{
(15) & \ds \sum_{n=1}^{+\infty} d_nF(nt) & = & \ds \frac{1}{2i\pi}\int_{\Re(s)=\frac{3}{2}} \zeta(s)^2\Gamma(s) F^{(-s)}(0)(-t)^{-s}ds\hfill\cr
 & & = & \ds \frac{1}{2i\pi}\int_{\Re(s)=-\frac{3}{2}} \zeta(-s)^2\Gamma(-s)F^{(s)}(0)(-t)^{s}ds\hfill\cr
(16) & & = & \ds -\frac{1}{t}\int_0^{+\infty}F(-u)(\ln(\frac{u}{-t})+2\gamma)du + \frac{F(0)}{4} + \sum_{n=1}^{N} \frac{\zeta(-2n+1)^2}{(2n-1)!}t^{2n-1} F^{(2n-1)}(0)\hfill\cr
 & & & {\hskip 4cm} + \ds \frac{1}{2i\pi}\int_{\Re(s)=2N} \Gamma(-s) \zeta(-s)^2F^{(s)}(0)(-t)^{s}ds.\hfill\cr
 }$$
 
The last equality $(16)$ is again justified by the integrals on the limit horizontal segments $[-\frac {3} {2} \pm iT, 2N \pm iT]$ being zero when $T\to +\infty$, for the reason explained between equation $(7)$ to $(8)$ (choose $n = 4N+2$).
\ms

This summation formula is to Voronoï's what the Euler-MacLaurin's formula is to Poisson's. For this reason, it is legitimate to call it the Euler-Voronoï's formula.
\bs

We move on towards Voronoï's formula. The approach is identical to that adopted for the Poisson formula. Let  $F\in \SM (]-\infty, 0], \R)$ be a function that we extend over $\R$ as an even function, and  $y > 0$.
\ms

Given the above, we can write

$$\matrix{
 & \ds H_1(y)   & = & \ds \sum_{n=1}^{+\infty} d_nF(ny) - \pa{\frac{1}{y}\int_0^{+\infty}F(u)(\ln(\frac{u}{y})
+2\gamma)du + \frac{F(0)}{4}}\hfill\cr
 & & = & \ds\frac{1}{2i\pi}\int_{\Re(s)=\frac{1}{2}} \Gamma(-s) \zeta(-s)^2F^{(s)}(0)y^{s}ds.   \hfill\cr
 }$$

Equations $ (11) $ and $ (12) $ give

$$\matrix{
\ds H_1(y) & = &\ds \frac{1}{2i\pi}\int_{\Re(s)=\frac{3}{2}}\Gamma(1-s) \zeta(1-s)^2F^{(s-1)}(0)y^{s-1}ds\hfill\cr
& = &\ds \frac{1}{2i\pi}\int_{\Re(s)=\frac{3}{2}}\zeta(s)^2\pa{2(2\pi)^{-2s}\Gamma(s)^2 \cos(\pi s/2)\zeta(s)^2\widehat{F}^{(-s)}(0)y^{s-1}}ds\hfill\cr}$$

From relation $(15)$ we get
$\ds H_1(y) = \frac{1}{y}\sum_{n=1}^{+\infty} d_nS\pa{\frac{n}{y}}$, with
$$\ds S(u) = \frac{1}{2i\pi}\int_{\Re(s)=\frac{3}{2}}2(2\pi)^{-2s}\Gamma(s)^2 \cos(\pi s/2)\widehat{F}^{(-s)}(0)u^{-s}ds.$$
\ms

Yet, if $L(y) = \ds \int_0^{+\infty} M\pa{\frac{y}{x}}N(x)\frac{dx}{x}$, then ${\cal M}[L](s) = {\cal M}[M](s){\cal M}[N](s)$ can be inverted into
$$L(x) = {\cal M}^{(-1)}[\Gamma(s)^2 M^{(-s)}(0)N^{(-s)}(0)](x). \leqno{(17)}$$
And in this case,
$$S(u) = \int_0^{+\infty} G\pa{\frac{u}{x}}\widehat{F}(x)\frac{dx}{x}$$
where $G$ verifies $\ds G^{(-s)}(0) = 2(2\pi)^{-2s}\cos\pa{s\pi/2}$. Although it is not in $\SM ([-\infty, a], \C) $, the solution $G (x) = 2\cos (4\pi^2x) $ is obvious: one just needs to be able to differentiate the cosine function!
\ms

Accordingly,
$\ds H_1(y) = \frac{2}{y}\sum_{n=1}^{+\infty} d_n\int_0^{+\infty}\cos\pa{\frac{4\pi^2 n}{xy}}\widehat{F}(x)\frac{dx}{x}$.
And finally,

$$\sum_{n=1}^{+\infty} d_nF(ny) = \frac{1}{y}\int_0^{+\infty}F(u)\pa{\ln\pa{\frac{u}{y}}+2\gamma}du + \frac{F(0)}{4} + \frac{2}{y}\sum_{n=1}^{+\infty} d_n\int_0^{+\infty}\cos\pa{\frac{4\pi^2 n}{xy}}\widehat{F}(x)\frac{dx}{x}.$$
\ms

It is the Voronoï's formula as soon as it is known that, for $\alpha >0$, the Fourier transform of $\ds \abs{x}^{-1}\cos(\alpha x^{-1})$ is $2K_0(2\sqrt{\alpha x}) - \pi Y_0(2\sqrt{\alpha x})$\, as a distribution.
\ms

If only from an aesthetic point of view, this writing of the Voronoï's formula is nicer than the original, but does not appear to be listed in the literature. Nevertheless, some characteristics of $F$ (especially his support) are less obvious here.
\bs
\ms


{\textsc{IX - Euler-Circle and Circle formulae}}
\ms

For an even function $F$ and under some hypotheses, we have the expression of the Circle [Mill]
$$F(0) + \sum_{n=1}^{+\infty} r_n F(n) = \pi\int_0^{+\infty}F(x)dx + \pi\pa{\sum_{n=1}^{+\infty} r_n\int_0^{+\infty}J_0(2\pi\sqrt{nx})F(x)dx} \leqno{(18)}$$
where \, $r_n = \# \{(a,b)\in\Z^2 \ / \ a^2+b^2=n\}$, $n\in\N$, and $J_0$ the first Bessel function.
\ms

Arithmetical reflections [Mill] lead to
$$\zeta(s)L(s,\chi_4) = \frac{1}{4}\sum_{n=1}^{+\infty} r_nn^{-s}, \ \ {\rm avec} \ \ L(s,\chi_4) = \sum_{k=0}^{+\infty} \frac{(-1)^k}{(2k+1)^s}.$$
Just like for $\zeta$, we have for this Dirichlet function a functional equation [Mill]:

$$\forall s\in\C, \ \ \ \ \xi_4(1-s) = \xi_4(s), \ \ {\rm with} \ \ \xi_4(s) = 2^s\pi^{-(s+1)/2}\Gamma(\frac{s+1}{2})L(s,\chi_4).\leqno{(19)}$$

To obtain $(18)$ from a contour integral, for $F\in \SM (]-\infty, 0], \C) $ we define,  for $t < 0$,
$$\matrix{
 & H(t) & = & \ds \frac{1}{2i\pi}\int_{\Re(s)=\frac{3}{2}} \zeta(s)L(s,\chi_4)\Gamma(s)F^{(-s)}(0)(-t)^{-s}ds\hfill\cr
(20) & & = & \ds \frac{1}{4}\sum_{n=1}^{+\infty} r_nF(nt)\hfill\cr
 & & = & \ds - \frac{L(1,\chi_4)}{t}F^{(-1)}(0)-\frac{L(0,\chi_4)}{2}F(0) + \frac{1}{2i\pi}\int_{\Re(s)=-\frac{1}{2}} \Gamma(s) \zeta(s)L(s,\chi_4)F^{(-s)}(0)(-t)^{-s}ds\hfill\cr
(21) & & = & \ds - \frac{\pi}{4t}F^{(-1)}(0)-\frac{1}{4}F(0) + \frac{1}{2i\pi}\int_{\Re(s)=-\frac{1}{2}} \Gamma(s) \zeta(s)L(s,\chi_4)F^{(-s)}(0)(-t)^{-s}ds.\hfill\cr}$$
Equation $(19)$ indeed shows that $L (-2n+1, \chi_4) = 0$ for any $n\in\N ^\star$. 
Since  $\zeta (-2n) = 0$, 
we conclude on the absence of pole for $\Gamma (s) \zeta (s) L (s, \chi_4) $ on $\Re(s) < 0$.
\ms

Consequently, for any $c > 0$, we have the functional equations  of $\xi$ and $\xi_4$
$$\matrix{
 & \ds \frac{1}{4}\sum_{n=1}^{+\infty} r_nF(nt) & = &\ds - \frac{\pi}{4t}F^{(-1)}(0) -\frac{1}{4}F(0) + \frac{1}{2i\pi}\int_{\Re(s)=c} \Gamma(-s) \zeta(-s)L(-s,\chi_4)F^{(s)}(0)(-t)^{s}ds\hfill\cr
(22) &  & = &\ds - \frac{\pi}{4t}F^{(-1)}(0) -\frac{1}{4}F(0) + \frac{1}{4i\pi}\int_{\Re(s)=c} \xi(s+1)\xi_4(s+1)F^{(s)}(0)(-t)^{s}ds\hfill\cr}$$

The regular part of formula's asymptotic expansion, which we can to call Euler-Circle, is therefore stationary from the order $0$.
\ms

Let us note that equalities $(20)$, $(21)$ and $(22)$ are again justified by the integrals on the limit horizontal segments  $[a \pm iT, b \pm iT]$ being zero, thanks to a similar argument as the one detailed in section V.
\bs

We continue our walk towards the Circle's formula with $y > 0$ and $F$ extended to $\R$ as an even function. For $y = -t >0$. From $(21)$ we deduct from $(21)$
$$H(y)- \frac{\pi}{4y}F^{(-1)}(0) + \frac{1}{4}F(0) = \ds \frac{1}{2i\pi}\int_{\Re(s)=\frac{3}{2}} \Gamma(1-s) \zeta(1-s)L(1-s,\chi_4)F^{(-1+s)}(0)y^{s-1}ds.$$

Then, just like for the Poisson's formula,
$$\matrix{
\ds H(y)- \frac{\pi}{4y}F^{(-1)}(0) + \frac{1}{4}F(0) & = & \ds \frac{1}{2i\pi}\int_{\Re(s)=\frac{3}{2}} \frac{\widehat{F}^{(-s)}(0)}{2\cos(\pi s/2)} 2\cos(\pi s/2)(2\pi)^{-s}\Gamma(s)\zeta(s)L(1-s,\chi_4)y^{s-1}ds.\hfill\cr
& = &\ds \frac{1}{2i\pi}\int_{\Re(s)=\frac{3}{2}} \Gamma(s)\widehat{F}^{(-s)}(0)(2\pi)^{-s}2^{2s-1}\frac{\pi^{-(s+1)/2}\Gamma(\frac{s+1}{2})}{\pi^{-(2-s)/2}\Gamma(\frac{2-s}{2})}\zeta(s)L(s,\chi_4)y^{s-1}ds.\hfill\cr
& = &\ds \frac{1}{2i\pi}\int_{\Re(s)=\frac{3}{2}} \frac{\sqrt{\pi}}{2}\pa{\frac{2}{\pi^2}}^{s}\Gamma(s)\widehat{F}^{(-s)}(0)\frac{\Gamma(\frac{s+1}{2})}{\Gamma(\frac{2-s}{2})}\zeta(s)L(s,\chi_4)y^{s-1}ds.\hfill\cr
}$$

Using $(20)$ this leads to $\ds H(y)- \frac{\pi}{4y}F^{(-1)}(0) + \frac{1}{4}F(0) = \frac{1}{4y}\pa{\sum_{n=1}^{+\infty} r_nS\pa{\frac{n}{y}}}$ with this time
$$S(u) = \frac{1}{2i\pi}\int_{\Re(s)=\frac{3}{2}} \Gamma(s)^2 \frac{\sqrt{\pi}}{2}\pa{\frac{2}{\pi^2}}^{s}\widehat{F}^{(-s)}(0)\frac{\Gamma(\frac{s+1}{2})}{\Gamma(s)\Gamma(\frac{2-s}{2})}u^{-s}ds.$$

The argument $(17)$ detailled for Voronoï's formula gives here
$$S(u) = \int_0^{+\infty}\widehat{F}(x) G\pa{\frac{u}{x}}\frac{dx}{x}$$
where $G$ statisfies 
$\ds G^{(-s)}(0) = \frac{\sqrt{\pi}}{2}\pa{\frac{2}{\pi^2}}^{s}\frac{\Gamma(\frac{s+1}{2})}{\Gamma(s)\Gamma(\frac{2-s}{2})} = \pi^{-2s}\sin\pa{\frac{\pi s}{2}}$
from Legendre's duplication formula.

Thus, $G(x) = -\sin(\pi^2 x)$ et $S(u) = \ds -\int_0^{+\infty}\widehat{F}(x) \sin\pa{\frac{\pi^2 u}{x}}\frac{dx}{x}$, and finally,
$$\sum_{n=0}^{+\infty} r_n F(ny) = \frac{\pi}{y}\int_0^{+\infty}F(u)du + \frac{1}{y}\sum_{n=1}^{+\infty} r_n\int_0^{+\infty}\sin\pa{\frac{\pi^2 n}{yu}}\widehat{F}(u)\frac{du}{u}. \leqno{(23)}$$

It is well known that, for $\alpha > 0$, $\pi J_0(2\sqrt{\alpha x})$ is the Fourier tranform of $x^{-1} \sin\pa{\alpha x^{-1}}$.
Hence, the Circle's formula can be derived from $(23)$ using Plancherel's formula. 
The last remark about the Voronoï's formula is still valid here: the equation $(23)$ is more aesthetic than Voronoï's formula 
and makes use of easier functions to be control (sinus versus $J_0$), 
but once more information on the support of $F$ is harder to see.
\ms
\bs


{\textsc{X - Euler-M\"obius-Poisson and M\"obius-Poisson's formulae}}
\ms

A formal M\"obius inversion of Poisson's formula would give
$$\sum_{n=1}^{+\infty} \frac{\mu_n}{n}F\pa{\frac{2\pi t}{n}} \approx -\frac{1}{2\pi t} \sum_{n=1}^{+\infty} \frac{\mu_n}{n} \widehat{F}\pa{\frac{1}{nt}}. \leqno{(24)}$$

Numerical tests on function $F (t) = e ^ {-t^2 / 2} $ reveal that this equality is true at first sight, but turns out to be wrong when precision increases.
\ms

A formal development of this inversion reads
$$2\pi\sum_{k=1}^{+\infty} \frac{F^{(2k)}(0)}{(2k)!} \frac{(2\pi)^{2k}t^{2k}}{\zeta(2k+1)} 
\approx 
\sum_{k=1}^{+\infty} \frac{\widehat{F}^{(2k)}(0)}{(2k)!} \frac{t^{-2k-1}}{\zeta(2k+1)}. \leqno{(25)}$$

The right-hand side together with previous sections give us
$$\frac{1}{2i\pi} \int_{C_\infty} \Gamma(-s)F^{(s)}(0) \frac{(2\pi)^{s}}{\zeta(s+1)} (-t)^{s}ds = \frac{1}{2i\pi} \int_{C_\infty} \Gamma(s)F^{(-s)}(0) \frac{(2\pi)^{-s}}{\zeta(-s+1)} (-t)^{-s}ds = 0$$
where $C_{\infty}$ is the limit of a sequence of rectangular paths.
\bs

To derive this equality, we assume that $F$ is even and is in $\SM(]-\infty,0],\C)$.
\ms

Let us give a name to the following hypotheses:
\be
\item[-]
$({\cal H}_0)$ There exist two reals $A>0$, $\alpha < \ds \frac{\pi}{2}$ satisfying for $s = \sigma + i\tau$, $\tau\in [-\frac{1}{2},\frac{3}{2}]$, $\abs{F^{(-s)}(0)} \ie Ae^{\alpha \abs{\tau}}$.
\item[-]
${\cal (H)}$ There exist two reals $A>0$, $\alpha < \ds \frac{\pi}{2}$ satisfying for $s = \sigma + i\tau$, $\tau\in \R$, \, $\abs{F^{(-s)}(0)} \ie Ae^{\alpha \abs{\tau}}$.
\ee

Thanks to Weil's formula [Weil], we prove the following lemma.
\ms

\textit{Lemma:} Let $F$ be a function satisfying ${\cal H}_0$.
Then, there exists a sequence of paths $(T_N)_{N\in\N}$, with $N \ie T_N\ie N+1$, such that for any $\theta \in\C$ 
such that $\alpha + \abs{\Re(\theta)} < \pi/2$, we have
$$\lim_{N\to +\infty}\int_{S_N} \abs{\frac{F^{(-s)}(0) e^{is\theta}}{2\cos(\pi s/2)\zeta(s)} }ds  = 0$$
where $(S_N)$ is one of the two sequences of segments $([-1/2 + iT_N,3/2 + iT_N])_{n\in\N}$ and $([-1/2 - iT_N,3/2 - iT_N])_{n\in\N}$.
\bs

\begin{proof}[Proof]
A delicate and classic approximation [Tit2] will prove to be greatly useful: if $s = \sigma + i\tau$ and if we note 
$\rho = \beta + i\gamma$ the non trivial generic zeroes of $\zeta$, then
$$\forall \sigma\in [-1,2], \ \ \ln(\zeta(s)) = \sum_{\abs{\tau-\gamma} \ie 1} \ln(s-\rho) + O(\ln(\tau)),$$
uniformly in $\sigma$, $\ln$ being the principal value of logarithm.

For $\ds \sigma\in[-\frac{1}{2},\frac{3}{2}]$ and $\abs{\tau} \se 1$, we can find $a > 0$ such that 
$$\ds \abs{\frac{1}{\zeta(s)}} \ie \abs{\tau}^a\abs{\prod_{\abs{\tau-\gamma} \ie 1} (s-\rho)}^{-1}.$$
Moreover, it is known [Tit2] that the number of zeroes in the strip $\Im(s) \in \pm[N,N+1]$ is also of order $O(\ln(N))$. Let $b > 0$ and $N_0$ be such that this number is upper-bounded by $b\ln(N)$ for $N \se N_0$.

For $N \se N_0$, we can find $T_N\in [N,N+1]$ such that $\abs{T_N - \gamma} \se \ds \frac{1}{2b\ln(N)}$ for all imaginary part $\gamma$ of a zero of $\zeta$.
\bs

Thus, for $s\in [-1 \pm iT_N,2 \pm iT_N]$, $\abs{\zeta(s)^{-1}} \ie (N+1)^a (2b\ln(N))^{b\ln(N)}$ and
$$\forall \sigma\in\left[-\frac{1}{2},\frac{3}{2}\right], \ \forall \theta \in\C, \ \ \abs{\frac{F^{(-s)}(0) e^{is\theta}}{2\cos(\pi s/2)\zeta(s)} } \ie \tilde{A} e^{-N\pi/2 + (N+1)(\alpha+\abs{\Re(\theta)})}(N+1)^a (2b\ln(N))^{b\ln(N)}$$
for some number $\tilde{A} > 0$,  which ends the lemma's proof. 
\end{proof}

For $\theta$ satisfying $\ds \abs{\Re(\theta)} < -\alpha + \frac{\pi}{2}$, let us define
$$L(s,\theta) = \frac{1}{2i\pi} \frac{1}{2\cos(\pi s/2)\zeta(s)} F^{(-s)}(0) e^{is\theta} = \frac{1}{2i\pi} \frac{(2\pi)^{-s}}{\zeta(-s+1)}\Gamma(s)F^{(-s)}(0) e^{is\theta}.$$
With the notations of  part V, if $C_{0,N} = R_{-\frac{1}{2},\frac{3}{2},T_N}$, we have

$$\lim_{N\to +\infty} \int_{C_{0,N}} L(s,\theta)ds = \lim_{N\to +\infty} \pa{ \int_{\frac{3}{2}-iT_N}^{\frac{3}{2}+iT_N} L(s,\theta)ds - \int_{-\frac{1}{2}-iT_N}^{-\frac{1}{2}+iT_N} L(s,\theta)ds }$$
which we call $\ds \int_{C_0,\infty} L(s,\theta)ds$ for brevity.
\bs

To get an asymptotic expansion, we want to push away the integration lines $\ds \Re(s) = {-1}/{2}$ and $\ds \Re(s) = {3}/{2}$. 
For this, let us notice that the previous lemma can similarly be applied to segments $([a\pm iT_N,b\pm iT_N])_{n\in\N}$, under hypothesis $ {\cal H} $, and this for any real numbers $a < b$.

On segments which real parts are included in $[3/2, + \infty [$, the fact that $\abs {\zeta (s)} ^ {-1}$ is bounded on $\Re (s) \se 3/2$ allows us to conclude. Similarly, the  $\zeta$ functional equation and  $\Gamma (1+i\tau) $ being decreasing allows one to prove that $\abs {\zeta (s)} ^ {-1}$ is bounded on $\Re (s) \ie-1/2$ and we can also conclude on segments which real parts are included in $]-\infty,-1/2] $. 
For any $n\in \N ^*$, an simple calculation of residue leads to

$$\lim_{N\to +\infty}\int_{-1/2-iT_N}^{-1/2+iT_N} L(s,\theta)ds = \sum_{k=1}^{n} \frac{F^{(2k)}(0)}{(2k)!} \frac{(2\pi)^{2k}e^{-2ik\theta}}{\zeta(2k+1)} + \lim_{N\to +\infty}\int_{-2n-1-iT_N}^{-2n-1+iT_N} L(s,\theta)ds$$
and
$$\lim_{N\to +\infty}\int_{3/2-iT_N}^{3/2+iT_N} L(s,\theta)ds = \frac{1}{2\pi } \sum_{k=1}^{n} \frac{\widehat{F}^{(2k)}(0)}{(2k)!} \frac{e^{(2k+1)i\theta}}{\zeta(2k+1)} + \lim_{N\to +\infty}\int_{2n+2-iT_N}^{2n+2+iT_N} L(s,\theta)ds.$$

The equality uses equation $(11)$ for even functions, and we obtain

$$\matrix{
\ds \int_{C_0,\infty} L(s,\theta)ds & = &\ds \frac{1}{2\pi } \sum_{k=1}^{n} \frac{\widehat{F}^{(2k)}(0)}{(2k)!} \frac{e^{(2k+1)i\theta}}{\zeta(2k+1)}-\sum_{k=1}^{n} \frac{F^{(2k)}(0)}{(2k)!} \frac{(2\pi)^{2k}e^{-2ik\theta}}{\zeta(2k+1)}\hfill\cr
& & {\hskip 1cm} \ds + \lim_{N\to +\infty}\int_{2n+2-iT_N}^{2n+2+iT_N} L(s,\theta)ds - \lim_{N\to +\infty}\int_{-2n-1-iT_N}^{-2n-1+iT_N} L(s,\theta)ds.\hfill\cr}$$

Finally, under rather broad assumptions, we note that $\ds \lim_{n\to +\infty}\int_{\pm(2n+1/2)-i\infty}^{\pm(2n+1/2)+i\infty} L(s,\theta)ds = 0$.

Thus, for $\ds \abs{\Re(\theta)} < \frac{\pi}{2} - \alpha$,
$$\frac{1}{2i\pi} \int_{C_0,\infty} \frac{F^{(-s)}(0)e^{is\theta}}{2\cos(\pi s/2)\zeta(s)}ds  = \ds \frac{1}{2\pi } \sum_{k=1}^{+\infty} \frac{\widehat{F}^{(2k)}(0)}{(2k)!} \frac{e^{(2k+1)i\theta}}{\zeta(2k+1)}-\sum_{k=1}^{+\infty} \frac{F^{(2k)}(0)}{(2k)!} \frac{(2\pi)^{2k}e^{-2ik\theta}}{\zeta(2k+1)}. \leqno{(26)}$$
\ms
which we can call the Euler-M\"obius-Poisson's formula.
\ms

Getting back to a positive real variable $y = e ^ {-i\theta} $, the Dirichlet development of $\zeta ^ {-1} $ and a lawful application of 
the Fubini theorem lead to
$$\forall y > 0, \ \ \frac{1}{2i\pi} \int_{C_0,\infty} \frac{F^{(-s)}(0)y^{-s}}{2\cos(\pi s/2)\zeta(s)}ds = \ds \frac{1}{2\pi y} \sum_{n=1}^{+\infty} \frac{\mu_n}{n} \widehat{F}\pa{\frac{1}{ny}} - \sum_{n=1}^{+\infty} \frac{\mu_n}{n}F\pa{\frac{2\pi y}{n}} \leqno{(27)}$$

\ms
which we can name the M\"obius-Poisson's formula this time. Let us remark that this latter equation uses the well-known equality presumably due to Euler $ \sum_{n}^{} {\mu_n}/{n} = 0$.  
\ms

As a final remark, let us point out that equations $(26)$ and $(27)$ highlight why formally inverting equations $(25)$ and $(24)$ is wrong.\bs

{\textsc{Prospective: under Riemann hypothesis}}

\bs

We now assume that we are working under hypotheses such that equations $(26)$ and $(27)$ are valid.
By assuming the Riemann hypothesis, stating that the non trivial zeroes of $\zeta$ are simple and located on the line 
$\Re(s) = {1}/{2} $, we note $Z$ the set of zeroes of $\zeta$.\ms

For $y > 0$, $(27)$ becomes

$$\sum_{\rho\in Z} \frac{F^{(-\rho)}(0)}{2\cos(\pi\rho/2)\zeta'(\rho)} y^{-\rho} = \ds \frac{1}{2\pi y} \sum_{n=1}^{+\infty} \frac{\mu_n}{n} \widehat{F}\pa{\frac{1}{ny}} - \sum_{n=1}^{+\infty} \frac{\mu_n}{n}F\pa{\frac{2\pi y}{n}}$$
which we make symmetrical using the change of variable $z = \sqrt{2\pi}y > 0$,
$$\sum_{\rho\in Z} \frac{F^{(-\rho)}(0)(2\pi)^{\rho/2}}{2\cos(\pi\rho/2)\zeta'(\rho)} z^{-\rho+1/2} = \frac{1}{\sqrt{2\pi z}} \sum_{n=1}^{+\infty} \frac{\mu_n}{n} \widehat{F}\pa{\frac{\sqrt{2\pi}}{n}\frac{1}{z}} - \sqrt{z}\sum_{n=1}^{+\infty} \frac{\mu_n}{n}F\pa{\frac{\sqrt{2\pi}}{n}z}. \leqno{(28)}$$

And for $\ds \abs{\Re(\theta)} < \frac{\pi}{2} - \alpha$, equation $(26)$ becomes
$$\sum_{\rho\in Z} \frac{F^{(-\rho)}(0)}{2\cos(\pi\rho/2)\zeta'(\rho)} e^{i\rho\theta} = \frac{1}{2\pi } \sum_{k=1}^{+\infty} \frac{\widehat{F}^{(2k)}(0)}{(2k)!} \frac{e^{i(2k+1)\theta}}{\zeta(2k+1)} - \sum_{k=1}^{+\infty} \frac{F^{(2k)}(0)}{(2k)!} \frac{(2\pi)^{2k}e^{-i2k\theta}}{\zeta(2k+1)}.$$

Again we symmetrise it by changing $\theta$ into $\theta -\frac{i}{2}\ln(2\pi)$:

$$\sum_{\rho\in Z} \frac{F^{(-\rho)}(0)(2\pi)^{\rho/2}}{2\cos(\pi\rho/2)\zeta'(\rho)} e^{i(\rho-\frac{1}{2})\theta} = \sum_{k=1}^{+\infty} \frac{(2\pi)^{k}}{(2k)!\zeta(2k+1)}\pa{\frac{1}{\sqrt{2\pi}}\widehat{F}^{(2k)}(0)e^{i(2k+\frac{1}{2})\theta} - F^{(2k)}(0)e^{-i(2k+\frac{1}{2})\theta}}. \leqno{(29)}$$

Let us note that sums over $Z$ are defined as $\ds \sum_{\rho\in Z} = \lim_{N\to +\infty} \sum_{\rho\in Z, \ \abs{\Im(\rho)}\ie T_N}$.
\bs

{\textsc{Applying equation $(29)$ to the test function $F(t) = e^{-t^2/2}$}}
\bs

One gets $\widehat{F}(t) = \sqrt{2\pi} F(t)$,\,  $\ds \frac{F^{(2n)}}{(2n)!}(0) = \frac{(-1)^n}{2^n n!}$\,  
and
$F^{(-s)}(0) = \ds \frac{\Gamma(s/2)}{\Gamma(s)} 2^{s/2 - 1}$.
\ms

We remind the following equalities
$$\forall N\in \N, \ \ \abs{\Gamma(N+1/2+i\tau)}^2 = \frac{\pi}{{\rm ch}(\pi \tau)}\abs{\prod_{k=0}^{N-1}((k+1/2)^2+\tau^2)}^2, $$
$$\forall N\in \N, \ \ \abs{\Gamma(-N-1/2+i\tau)}^2 = \frac{\pi}{{\rm ch}(\pi \tau)}\abs{\prod_{k=0}^{N}((k+1/2)^2+\tau^2)}^{-2}.$$
\ms

The validity of $(29)$ in that case is verified through the following steps ($s =\sigma+i\tau$):
\ms

\be
\item[$\bullet$]
$F$ verifies hypothesis ${\cal H}$, and for any $\alpha < \ds \frac{\pi}{4}$, we can find $A > 0$ such that $\abs{F^{(-s)}(0)} \ie A2^{\abs{\sigma}/2}e^{\alpha \abs{\tau}}$.
\ms

\item[$\bullet$]
For any $\abs{\Re(\theta)} < \ds \frac{\pi}{4}$, \,
$\ds \lim_{n\to +\infty}\int_{2n-i\infty}^{2n+i\infty} \frac{2^{s/2 - 1}}{2\cos(\pi s/2)\zeta(s)} \frac{\Gamma(s/2)}{\Gamma(s)} e^{is\theta} ds = 0$.
\ms

Since $s = 2n+i\tau$, for $\theta = \theta_1+i\theta_2$ 
with $\abs{\theta_1} < \pi/4$, 
we have
\ms
$$\abs{\frac{2^{s/2 - 1}}{2\cos(\pi s/2)\zeta(s)} \frac{\Gamma(s/2)}{\Gamma(s)} e^{is\theta}} = \abs{\frac{\pi^{1/2} 2^{-s/2} e^{-\tau\theta_1 - 2n\theta_2}}{2\cos(\pi s/2)\zeta(s)\Gamma((s+1)/2)}} \ie \frac{\pi^{1/2} 2^{-n} e^{-\tau\theta_1 - 2n\theta_2}}{{\rm ch}(\pi \tau/2)\abs{\Gamma((n+1/2+i\tau/2)}}$$

and the well-known equality (reflection formula and functional equation of $\Gamma$)
$$\forall n\in \N^*, \ \forall \tau\in\R, \ \ \abs{\Gamma(n+1/2+i\tau/2)}^2 = \frac{\pi}{{\rm ch}(\pi \tau/2)}\abs{\prod_{k=0}^{n-1}((k+1/2)^2+(\tau/2)^2)}^2$$
allows us to easily conclude.
\ms

\item[$\bullet$]
For any $\abs{\Re(\theta)} < \ds \frac{\pi}{4}$, \,
$\ds \lim_{n\to +\infty}\int_{-(2n+1)-i\infty}^{-(2n+1)+i\infty} \frac{(2\pi)^{-s}2^{s/2 - 1}}{\zeta(-s+1)} \Gamma(s/2) e^{is\theta} ds = 0$.
\ms

The argument is similar to the previous one, but this time using the equality
$$\forall n\in \N^*, \ \forall \tau\in\R, \ \ \abs{\Gamma(-n-1/2+i\tau/2)}^2 = \frac{\pi}{{\rm ch}(\pi \tau/2)}\abs{\prod_{k=0}^{n}((k+1/2)^2+(\tau/2)^2)}^{-2}.$$
\ee
\bs

Defining $\rho = \frac{1}{2} + i\tau$ for $\rho\in Z$ ($\tau$ is therefore assumed to be real), equality $(29)$ becomes
$$\frac{i}{2}\sum_{\rho\in Z} \frac{2^{\rho}\pi^{\rho/2}\Gamma(\rho/2)}{4\cos(\pi\rho/2)\Gamma(\rho)\zeta'(\rho)} e^{\tau\theta} = \sum_{k=1}^{+\infty} \frac{(-1)^k \pi^{k}}{k!\ \zeta(2k+1)}\sin((2k+1/2)\theta).$$
\ms

Let us call $Z^+$ the set of elements of $Z$ having positive imaginary part. We know that $Z = Z^+ \cup \{1-\rho, \ \rho\in Z^+\}$.
\ms

The $\zeta$ functional equation gives, for $\rho\in Z$: $\Gamma(\frac{\rho}{2})\zeta'(\rho) = - \pi^{\rho-\frac{1}{2}}\Gamma(\frac{1-\rho}{2})\zeta'(1-\rho)$. 
A direct calculation then gives, for $\rho\in\Z^+$, $C(\rho) = -C(1-\rho)$, where $\ds C(s) = i \frac{2^{\rho}\pi^{s/2}\Gamma(s/2)}{\cos(\pi s/2)\Gamma(s)\zeta'(s)}$. 
\ms

Finally, we easily prove that $C(\overline{z}) = -\overline{C(z)}$ and $C(\rho)\in\R$ for $\rho\in Z$ since $\overline{\rho} = 1-\rho$ by hypothesis.
\ms

Consequently,
$$i\sum_{\rho\in Z^+} \frac{2^{\rho}\pi^{\rho/2}\Gamma(\rho/2)}{4\cos(\pi\rho/2)\Gamma(\rho)\zeta'(\rho)} {\rm sh}(\tau\theta) = \sum_{k=1}^{+\infty} \frac{(-1)^k \pi^{k}}{k!\ \zeta(2k+1)}\sin((2k+1/2)\theta). \leqno{(30)}$$
\ms

Numerical simulations confirm the convergence of the left-hand side term when $\ds \abs{\Re{(\theta)}} < {\pi}/{4}$, 
the equality under the same hypothesis, and the divergence of the  the left-hand side for $\ds \abs{\Re(\theta)} > {\pi}/{4}$. 
\bs
\ms

The author expresses his gratitude to Alban Levy, the bilingual manuscript's proofreader.
\bs
\ms



{\textsc{Bibliography:}}
\ms

\bib{Edwa}{H.M. Edwards, {\it Riemann's Zeta Function}, Academic Press, New York, 1974.}

\bib{Endr}{S. Endres and F. Steiner, {\it A new proof of the Voronoï summation formula}, J. Phys. A: Math. Theor., 2011.}

\bib{Heja}{D. A. Hejhal., {\it A note on the Voronoï summation formula}, Monatshefte für Mathematik, 1979.}

\bib{Mill}{S. Miller and W. Schmid, {\it Summation formulae, from Poisson and Voronoi to the Present}, Noncommutative Analysis, in Progress in Mathematics 220, Birkhäuser, 2003.}

\bib{Ross}{K. S. Miller, B. Ross (Editor), {\it An Introduction to the Fractional Calculus and Fractional Differential Equations}, John Wiley and Sons, 1993.}

\bib{Tit1}{E. C. Titchmarsh, {\it Introduction to the Theory of Fourier Integrals}, Chelsea Publishing Compagny, 1986. Réédition de la seconde édition, Oxford, 1948.}

\bib{Tit2}{ E. C. Titchmarsh, {\it The Theory of the Riemann Zeta-Function} (revised by D. R. Heath-
Brown), Clarendon, Oxford, 1986.}

\bib{Weil}{A. Weil, {\it Sur les "formules explicites" de la théorie des nombres premiers}, Comm. Sém. Math. Univ. Lund, 1952.}

\end{document}